\documentclass[10pt,twoside]{article}
\usepackage{amsmath}
\usepackage{amsthm}
\usepackage{Latex-document}
\usepackage{amssymb}
\usepackage[mathscr]{eucal}

\makeatletter
\def\@seccntformat#1{\csname the#1\endcsname.\quad}
\renewcommand{\thesection}{\@arabic\c@section}
\renewcommand{\thesubsection}{\thesection.\arabic{subsection}}

\makeatother

\newtheorem{Theorem}[equation]{Theorem}

\theoremstyle{definition}

\theoremstyle{remark}

\newtheorem{Remarks}[equation]{Remarks}

\newcommand{\thmref}[1]{Theorem~\ref{#1}}
\newcommand{\secref}[1]{\S\ref{#1}}

\newcommand{\subsecref}[1]{\S\ref{#1}}

\newcommand{\defeq}{\overset{\operatorname{\scriptstyle def.}}{=}}
\newcommand{\C}{{\mathbf C}}
\newcommand{\Z}{{\mathbf Z}}

\newcommand{\R}{{\mathbf R}}
\newcommand{\proj}{{\mathbf P}}
\newcommand{\CP}{\proj}

\newcommand{\SL}{\operatorname{\rm SL}}

\newcommand{\GL}{\operatorname{GL}}

\newcommand{\g}{{\mathfrak g}}
\newcommand{\ag}{{\widehat{\mathfrak g}}}

\newcommand{\rank}{\operatorname{rank}}

\newcommand{\ve}{\varepsilon}
\newcommand{\M}{{\mathfrak M}} 

\newcommand{\La}{{\mathfrak L}} 
\newcommand{\bv}{{\mathbf v}} 
\newcommand{\bw}{{\mathbf w}} 
\newcommand{\bC}{{\mathbf C}} 

\newcommand{\codim}{\operatorname{codim}} 
\newcommand{\topdeg}{\operatorname{top}} 

\newcommand{\Ua}{{\mathbf U}_q(\widehat{\mathfrak g})} 
\newcommand{\Ul}{{\mathbf U}_q({\mathbf L}{\mathfrak g})} 

\newcommand{\Ule}{{\mathbf U}_{\varepsilon}({\mathbf L}\mathfrak g)}
\newcommand{\bU}{\mathbf U} 

\newcommand{\shfO}{\mathcal O}

\newcommand{\Wedge}{{\textstyle \bigwedge}}
\newcommand{\Zw}{Z(\bw)}

\newcommand{\Hilb}[2]{\operatorname{Hilb}^{#2}{#1}} 
\newcommand{\Hilbn}[1]{\Hilb{#1}{n}}
\newcommand{\Ut}{{\mathbf U}_q({\mathbf L}{\ag})} 
\newcommand{\Uti}{{\mathbf U}^{\Z}_q({\mathbf L}{\ag})}
\newcommand{\Ute}{{\mathbf U}_{\varepsilon}({\mathbf L}\ag)}

\def\volumeno{I}

\title{\bf Geometric Construction of \vskip -2mm Representations of Affine
Algebras\vskip 6mm}
\author{Hiraku Nakajima\vspace*{-0.5cm}\thanks{Department of Mathematics,
Kyoto University, Kyoto 606-8502, Japan. E-mail: nakajima@kusm.kyoto-u.ac.jp}}
\date{\vspace{-8mm}}

\begin{document}

\maketitle

\thispagestyle{first} \setcounter{page}{423}

\begin{abstract}\vskip 3mm

Let $\Gamma$ be a finite subgroup of $\SL_2(\C)$. We consider
$\Gamma$-fixed point sets in Hilbert schemes of points on the affine
plane $\C^2$. The direct sum of homology groups of components has a
structure of a representation of the affine Lie algebra $\ag$
corresponding to $\Gamma$. If we replace homology groups by
equivariant $K$-homology groups, we get a representation of the
quantum toroidal algebra $\Ut$. We also discuss a higher rank
generalization and character formulas in terms of intersection
homology groups.

\vskip 4.5mm

\noindent {\bf 2000 Mathematics Subject Classification:} 17B37,
14D21, 14L30, 16G20, 33D80.

\noindent {\bf Keywords and Phrases:} Affine Lie algebras, Quantum
toroidal algebras, Hilbert schemes, Quiver varieties.
\end{abstract}

\vskip 12mm

\section{Finite subgroups of {\boldmath $\SL_2(\C)$} and simple Lie
algebras}\label{sec:intor}
\setzero\vskip-5mm \hspace{5mm }

Let $\Gamma$ be a finite subgroup of $\SL_2(\C)$. The classification
of such subgroups has been well-known to us, since they are
essentially symmetry groups of regular polytopes. They are cyclic
groups, binary dihedral groups, and binary polyhedral groups (Klein (1884)).

It has been also known that we can associate a complex simple Lie
algebra $\g$ to $\Gamma$. This can be done in two ways. The first one
is geometric and due to DuVal (1934). The second one is algebraic, and
is due to McKay (1979).

Let us explain the two constructions and subsequent developments
briefly. More detailed account can be found in \cite{IN}.

\subsection{Minimal resolution of {\boldmath $\C^2/\Gamma$}}\label{subsec:DuVal}
\vskip-5mm \hspace{5mm }

Let us consider the quotient space $\C^2/\Gamma$. This space has an
isolated singularity at the origin. We have a unique {\it minimal\/}
resolution $\pi\colon M\to \C^2/\Gamma$, in the sense that all other
resolutions factor through $\pi$. (For general singularities, we have
many resolutions. This speciality occurs in $2$-dimensional case.)
This singularity is called a {\it simple singularity}, and has
been intensively studied from various points of view. In
particular, the following are known (see e.g., \cite{BPV}):
\begin{enumerate}
\item The exceptional set $\pi^{-1}(0)$ consists of the union of
projective lines.
\item We draw a diagram so that vertices correspond to projective
lines (irreducible components) and two vertices are connected by an
edge if they intersect. Then we obtain a Dynkin diagram of type $ADE$.
\end{enumerate}
We thus have bijections
\begin{equation*}
   \left\{\text{irreducible components of $\pi^{-1}(0)$}\right\}
   \longleftrightarrow
   \left\{\text{vertices of the Dynkin diagram}\right\}.
\end{equation*}

The Dynkin diagram appears in the classification of simple Lie
algebras. Thus we have a complex simple Lie algebra $\g$ corresponding
to $\Gamma$. Since vertices of the Dynkin diagram correspond to simple
coroots of $\g$, the above bijection gives an isomorphism (of vector
spaces)
\begin{equation}\label{eq:Cartan}
   \mathfrak h\xrightarrow{\sim} H_2(\pi^{-1}(0),\C),
\end{equation}
where $\mathfrak h$ is the complex Cartan subalgebra of $\g$.

This correspondence $\Gamma \to \g$ is based on the classification of
simple Lie algebras since they attach a Dynkin diagram to $\Gamma$. So
the reason why such a result holds remained misterious.
A deeper connection between two objects were conjectured by
Grothendieck, and obtained by Brieskorn (1970) and Slodowy
(1980). They constructed the simple singularity $\C^2/\Gamma$ in $\g$.
Moreover, its semi-universal deformation and a simultaneous resolution
were also constructed using geometry related to $\g$. We do not recall
their results here, so the interested reader should consult \cite{Sl}.

\subsection{McKay correspondence}\label{subsec:McKay}
\vskip-5mm \hspace{5mm }

Let $\{ \rho_i \}_{i\in I}$ be the set of (isomorphism classes of)
irreducible representations of $\Gamma$. It has a special element
$\rho_0$, the class of trivial representation. Let $Q$ be the
$2$-dimensional representation given by the inclusion $\Gamma\subset
\SL_2(\C)$.  Let us decompose $Q \otimes \rho_i$ into irreducibles, $Q
\otimes \rho_i = \bigoplus_j a_{ij} \rho_j$, where $a_{ij}$ is the
multiplicity. We draw a diagram so that vertices correspond to
$\rho_i$'s, and there are $a_{ij}$ edges between $\rho_i$ and
$\rho_j$. (Note that $a_{ij} = a_{ji}$ thanks to the self-duality of
$Q$).
Then McKay \cite{Mc} observed that the graph is an affine Dynkin
diagram of $\tilde{A}_n^{(1)}, \tilde{D}_n^{(1)}, \tilde{E}_6^{(1)},
\tilde{E}_7^{(1)}$ or $\tilde{E}_8^{(1)}$, i.e., the Dynkin diagram of
an untwisted affine Lie algebra $\ag$ attached to a simple Lie algebra
$\g$ of type $ADE$.
Furthermore it is also known that the Dynkin diagram given in the
previous subsection is obtained by the affine Dynkin diagram by removing
the vertex corresponding to the trivial representation $\rho_0$.
We thus have bijections
\begin{equation*}
   \left\{\text{irreducible representations of $\Gamma$}\right\}
   \longleftrightarrow
   \left\{\text{vertices of the affine Dynkin diagram}\right\}.
\end{equation*}

The original McKay's proof was based on the explicit calculation of
characters. The reason why such a result holds remained misterious
also in this case. A geometric explanation via the $K$-theory of the
minimal resolution $M$ of $\C^2/\Gamma$ was subsequently given by
Gonzalez-Sprinberg and Verdier \cite{GV}. In particular, they proved
that there exists a natural geometric construction of an isomorphism
(of abelian groups)
\begin{equation*}
   R(\Gamma) \xrightarrow{\sim} K(M),
\end{equation*}
where $R(\Gamma)$ is the representation ring of $\Gamma$, and $K(M)$
is the Grothendieck group of the abelian category of algebraic vector
bundles over $M$. This result is strengthened and generalized to
the higher dimensional case $\Gamma\subset\SL_3(\C)$ \cite{BKR}.

Note that the above isomorphism together with the Chern character
homomorphism leads to an isomorphism
\(
   R(\Gamma)\otimes_\Z\C \xrightarrow{\sim} H^*(M,\C),
\)
which gives an isomorphism
\begin{equation}\label{eq:affine}
   R(\Gamma)\otimes_\Z\C \xrightarrow{\sim}
   \left(\mathfrak h\oplus \C h_0\right)^*,
\end{equation}
combined with \eqref{eq:Cartan}. Here $h_0$ is the $0$th simple coroot
of the affine Lie algebra $\ag$, and corresponds to the dual of the
trivial representation $\rho_0$. It corresponds to $H_0(M,\C) \cong
H_0(\pi^{-1}(0),\C)$.

Compared with correspondence in \subsecref{subsec:DuVal}, our
situation is less satisfactory: we only get $\mathfrak h$ and the role
of $\g$ or $\ag$ is less clear. This is the starting point of our
whole construction. We construct $\ag$ entirely from $\Gamma$ in some
sense. For another approach, see \cite{FJW}.

\section{Hilbert schemes of points and their {\boldmath $\Gamma$}-fixed point
components --- quiver varieties}\label{sec:Hilb}
\setzero\vskip-5mm \hspace{5mm }

In 1986, Kronheimer \cite{Kr} constructed a simple singularity
$\C^2/\Gamma$, its deformation and simultaneous resolution, i.e.,
those spaces constructed by Brieskorn-Slodowy by a totally different
method. His construction is based on the theory of `quivers', which is
a subject in noncommutative algebras. (See also \cite{CS} for a
different approach.) Subsequently in 1989, Kronheimer and the author
\cite{KN} gave a description of moduli spaces of instantons (and
coherent sheaves) on those spaces in terms of a quiver. It is an
analog of the celebrated ADHM description of instantons on $S^4$. In
1994, this description was further generalized under the name of
`quiver varieties' by the author \cite{Na:1994}.
The purpose of this and next sections is to define quiver varieties
from a slightly different point of view. This is a most economical
approach to introduce quiver varieties, while it does not explain why
it is something to do with quivers.

Let $\Hilbn{(\C^2)}$ be the Hilbert scheme of $n$ points in the affine
plane $\C^2$. As a set, it consists of ideals $I$ of the polynomial
ring $\C[x,y]$ such that the quotient $\C[x,y]/I$ has dimension $n$ as
a vector space. Grothendieck constructed $\Hilbn{(\C^2)}$ as a
quasi-projective scheme (for more general setting), but we do not go
to this direction in detail. A typical point of $\Hilbn{(\C^2)}$ is an
ideal of functions vanishing at $n$ distinct points in $\C^2$. The
space parametrizing (unordered) $n$ distinct points is an open subset
of the $n$th symmetric product $S^n(\C^2) = (\C^2)^n/S_n$ of $\C^2$,
where $S_n$ is the symmetric group of $n$ letters acting on $(\C^2)^n$
by permutation of factors. The symmetric product parametrises
unordered $n$ points with multiplicities. The Hilbert scheme
$\Hilbn{(\C^2)}$ is a different completion of the open set.
Two completions are related: Mapping $I$ to its support counted with
multiplicities, we get a morphism $\pi\colon\Hilbn{(\C^2)}\to
S^n(\C^2)$ which is called a {\it Hilbert-Chow morphism}. We have
following important geometric results on $\Hilbn{(\C^2)}$:
\begin{enumerate}
\item $\Hilbn{(\C^2)}$ is a resolution of singularities of
$S^n(\C^2)$ (Fogarty).
\item $\Hilbn{(\C^2)}$ has a holomorphic symplectic structure
(Beauville, Mukai).
\end{enumerate}
In fact, the author constructed a hyper-K\"ahler structure on
$\Hilbn{(\C^2)}$, which induces Beauville-Mukai's symplectic form, by
describing it as a hyper-K\"ahler quotient.  See \cite{Lecture} and
G\"ottsche's article in this ICM proceeding for more recent results on
$\Hilbn{(\C^2)}$.

Let $\Gamma$ be a finite subgroup of $\SL_2(\C)$ as above. Its natural
action on $\C^2$ induces an action on $\Hilbn{(\C^2)}$ and $S^n(\C^2)$
such that the Hilbert-Chow morphism $\pi$ is $\Gamma$-equivariant. Let
us consider the fixed point set $\Hilbn{(\C^2)}^\Gamma$,
$(S^n(\C^2))^\Gamma$. The latter is easy to describe:
\begin{equation*}
  (S^n(\C^2))^\Gamma = S^m (\C^2/\Gamma),
\end{equation*}
where $m$ is the largest integer less than or equal to
$n/\#\Gamma$. The difference $n - m\#\Gamma$ is the multiplicity of the
origin. The former space $\Hilbn{(\C^2)}^\Gamma$ is a union of
nonsingular submanifolds of $\Hilbn{(\C^2)}$. If
$I\in\Hilbn{(\C^2)}^\Gamma$, the quotient $\C[x,y]/I$ has a structure
of a representation of $\Gamma$. For an isomorphism class $\bv$ of a
representation of $\Gamma$, we define $M(\bv)$ as
\begin{equation*}
   M(\bv) = \left\{\left. I \in \Hilbn{(\C^2)}^\Gamma\, \right|
    \left[\C[x,y]/I \right] = \bv \right\},
\end{equation*}
where $\left[\C[x,y]/I \right]$ is the isomorphism class of
$\C[x,y]/I$.  Since isomorphism classes are parametrized by discrete
data, i.e., dimensions of isotropic components, the isomorphism class
of $\left[\C[x,y]/I \right]$ is constant on each connected
component. Therefore $M(\bv)$ is a union of connected component. In
fact, Crawley-Boevey recently proves that $M(\bv)$ is connected (in
fact, he proved it for more general case including varieties discussed
in next section) \cite{CB}.
Moreover, $M(\bv)$ has induced holomorphic symplectic and
hyper-K\"ahelr structures. It is an example of quiver varities of
affine type. (See remark at the end of the next section.)

The simplest but nontrivial example is the case when $\bv$ is the
class of the regular representation of $\Gamma$.
Under \eqref{eq:affine}, the regular representation corresponds to the
imaginary root $\delta$, which is the positive generator of the kernel
of the affine Cartan matrix, is identified with the dimension vector
of the regular representation of $\Gamma$.
The dimension of the regular representation is equal to $\#\Gamma$,
and thus the fixed point set in the symmetric product is
$(S^{\#\Gamma}(\C^2))^\Gamma = \C^2/\Gamma$. We can consider this as
the space of $\Gamma$-orbits. A typical point is a free
$\Gamma$-orbit, and is also a point in $\Hilb{(\C^2)}{\#\Gamma}^\Gamma$ as
the ideal vanishing at the orbit. In fact, it is not difficult to see
that $M(\bv)$ is isomorphic to the minimal resolution $M$ of
$\C^2/\Gamma$.  The resolution map $\pi\colon M\to \C^2/\Gamma$ is
given by the restriction of the Hilbert-Chow morphism. This result was
obtained by Ginzburg-Kapranov (unpublished) and Ito-Nakamura \cite{IN}
independently, but is also a re-interpretation of Kronheimer's
construction \cite{Kr}. The precise explanation was given in
\cite[Chapter 4]{Lecture}.

Recently higher dimensional $M(\bv)$ attract attention in
connection with the McKay correspondence for wreath products
$\Gamma\wr S_n$ \cite{Wang,Kz,Ha}. These $M(\bv)$ are
diffeomorphic to the Hilbert schemes of points on the minimal
resolution $\Hilbn{M}$.

\section{A higher rank generalization of Hilbert schemes}
\setzero\vskip-5mm \hspace{5mm }
\newcommand{\linf}{\ell_\infty}

We give a higher rank generalization of Hilbert schemes in this
section. But geometric structures remain unchanged for general
cases. So a reader, who wants to catch only a {\it rough\/} picture,
could safely skip this section.

Let $\CP^2$ be the projective plane with a fixed line $\linf$. So
$\CP^2 = \C^2\sqcup \linf$. Let $\M(n,r)$ be the framed moduli space
of torsion free sheaves on $\CP^2$ with rank $r$ and $c_2=n$, i.e. the
set of isomorphism classes of pairs $(E,\varphi)$, where
$E$ is a torsion free sheaf of $\rank E =r$, $c_2(E)=n$, which is
locally free in a neighbourhood of $\linf$, and $\varphi$ is an
isomorphism
\(
  \varphi \colon E|_{\linf} \overset{\sim}{\to} {\cal O}_{\linf}^{\oplus r}
\)
(framing at infinity). It is known that this space has a structure of
a quasi-projective variety \cite{HL}. This is a higher rank
generalization of the Hilbert scheme $\Hilbn{(\C^2)}$. The analog of
$\C[x,y]/I$ is $H^1(\CP^2, E(-1))$ and it is known that
$H^0(\CP^2, E(-1)) = H^2(\CP^2, E(-1)) = 0$ \cite[Chapter 2]{Lecture}.
It is also known that $\M(n,r)$ has a holomorphic symplectic (in fact,
hyper-K\"ahler) structure \cite[Chapter 3]{Lecture}.

The higher rank generalization of the symmetric product $S^n(\C^2)$ is
the so-called Uhlenbeck compactification of the framed moduli space of
locally free sheaves. (On the other hand, $\M(n,r)$ is called
Gieseker-Maruyama compactification.) It is
\begin{equation*}
   \M_0(n,r) = \bigsqcup_{n'+n'' = n}
   \M_0^{\operatorname{reg}}(n',r) \times S^{n''} \C^2,
\end{equation*}
where $\M_0^{\operatorname{reg}}(n',r)$ is the open subset of
$\M(n',r)$ consisting of framed {\it locally free\/} sheaves
$(E,\varphi)$. It is known that $\M_0(n,r)$ has a structure of an
affine algebraic variety \cite[Chapter 3]{DK}. Moreover, the map
\begin{equation*}
   (E,\varphi) \mapsto
   \left(E^{\vee\vee}, \varphi,
     \operatorname{Supp}(E^{\vee\vee}/E)\right)
\end{equation*}
gives a projective morphism $\pi\colon \M(n,r)\to \M_0(n,r)$
\cite[Chapter 3]{Lecture}, where $E^{\vee\vee}$ is the double dual
of $E$, which is locally free on surfaces, and
$\operatorname{Supp}(E^{\vee\vee}/E)$ is the support of
$E^{\vee\vee}/E$, counted with multiplicities.

When $r = 1$, there exists only one locally free sheaf which is
trivial at $\linf$, i.e., the trivial line bundle ${\cal O}_{\CP^2}$. So
the first factor of the above disappears : $\M_0(n,1) = S^n
\C^2$. Moreover, for $E\in\M(n,1)$, the double dual $E^{\vee\vee}$
must be the trivial line bundle by the same reason. It means that $E$
is an ideal sheaf of the structure sheaf ${\cal O}_{\CP^2}$, so is a
point in the Hilbert scheme $\Hilbn{(\C^2)}$. Thus we recover the
situation studied in \secref{sec:Hilb}.

Let $\Gamma$ be a finite subgroup of $\SL_2(\C)$ as before. We take
and fix a lift of the $\Gamma$-action to ${\cal O}_{\linf}^{\oplus
r}$. It is written as $W\otimes_\C {\cal O}_{\linf}$, where $W$ is a
representation $W$ of $\Gamma$. We denote by $\bw$ the isomorphism
class of $W$ as before. Now $\Gamma$ acts on $\M(n,r)$, $\M_0(n,r)$
and we can consider the fixed point sets $\M(n,r)^\Gamma$,
$\M_0(n,r)^\Gamma$. We decompose the former as
\begin{equation*}
   \M(n,r)^\Gamma = \bigsqcup_\bv \M(\bv,\bw),
\end{equation*}
where $\M(\bv,\bw)$ consists of the framed torsion free sheaves
$(E,\varphi)$ such that the isomorphism class of $H^1(\CP^2,E(-1))$,
as a representation of $\Gamma$, is $\bv$. Each $\M(\bv,\bw)$, if it
is nonempty, inherits a holomorphic symplectic and hyper-K\"ahler
structures from $\M(n,r)$.

Arbitrary quiver varieties of affine types with complex parameter
equal to $0$ are some $\M(\bv,\bw)$. The identification with
original definition was implicitly given in \cite{Lecture}. It was
independently rediscovered by Lusztig \cite{Lu:quiver}. See also
\cite{VV:cyc}.
Arbitrary quiver varieties of affine types with {\it nonzero\/}
complex parameter are also important in representation theory
\cite{EG}, though we do not discuss here.
Original definition of the varieties was given in terms of quivers.
Later these were identified with framed moduli spaces of instantons on
a noncommutative deformation $\R^4$ \cite{NS} or those of torsion free
sheaves on a noncommutative deformation of $\CP^2$ \cite{KKO,BGK}.

\section{Stratification and fibers of {\boldmath $\pi$}}\label{sec:stratum}
\setzero\vskip-5mm \hspace{5mm }

This technical section will be used to state character formulas
later. A reader who only want to know only a {\it rough\/} picture
can be skip this section.

We have the following stratification of $(S^n(\C^2))^\Gamma$ and
its higher rank analog $\M_0(n,r)^\Gamma$. The space
$\M_0(n,r)^\Gamma$ also decompose as
\begin{equation}\label{eq:strat}
\begin{gathered}
   (S^n(\C^2))^\Gamma = \bigsqcup_{m\le n} S^m_\lambda (\C^2/\Gamma),
\\
   \M_0(n,r)^\Gamma = \bigsqcup_{\substack{\bv^0,\lambda \\
   m + |\bv^0| \le n}}
   \M_0^{\operatorname{reg}}(\bv^0,\bw) \times S^m_\lambda (\C^2/\Gamma),
\end{gathered}
\end{equation}
where $\M_0^{\operatorname{reg}}(\bv^0,\bw)$ is defined exactly as
above (it is possibly an empty set), $|\bv^0|$ is the dimension of
$\bv^0$ as a complex vector space, $\lambda =
(\lambda_1,\dots,\lambda_r)$ is a partition of $m$ and
\begin{equation*}
   S^m_\lambda (\C^2/\Gamma) =
   \left\{ \left.\sum_{i=1}^r \lambda_i [x_i] \in S^m (\C^2/\Gamma)\,
  \right|
  \text{$x_i\neq 0$ and $x_i \neq x_j$ for $i \neq j$} \right\}.
\end{equation*}
The differences $n-m$ and $n - (m+|\bv^0|)$ are the multiplicity of
the cycle at the origin $0$.

Now it becomes clear that the case $(S^n(\C^2))^\Gamma$ is the special
case of $\M_0(n,r)^\Gamma$ with $\bw = \rho_0$, $\bv^0 = 0$. So from
now, we only consider the second case.

For $x\in\M_0(n,r)^\Gamma$, let $\M(\bv,\bw)_x$ be the inverse image
$\pi^{-1}(x)$ in $\M(\bv, \bw)$.
The most important one is the central fiber, i.e., the fiber over
\begin{equation*}
   x = (W\otimes_\C \shfO_{\CP^2}, \varphi, 0).
\end{equation*}
In this case, we use the special notation $\La(\bv,\bw)$. It is
known that this is a Lagrangian subvariety of $\M(\bv,\bw)$.
Suppose that $x = (E_0,\varphi,C)$ is contained in the
stratum
\(
   \M_0^{\operatorname{reg}}(\bv^0,\bw)
     \times S^m_\lambda (\C^2/\Gamma)
\).
Then the fiber $\M(\bv,\bw)_x$ is a pure dimensional subvariety
in $\M(\bv,\bw)$, which is a product of
$\La(\bv_s,\bw_s)$ and copies of punctual Hilbert schemes
$\operatorname{Hilb}_0^{\lambda_i}{(\C^2)}$ for some $\bv_s$, $\bw_s$.
The proof of this statement in \cite[\S6]{Na:1994},
\cite[\S3]{Na-qaff} was given only when $m = 0$ and explained in terms
of quivers, so we give more direct argument in our situation.
The fiber $\M(\bv,\bw)_x$ parametrises $\Gamma$-invariant
subsheaves $E$ of $E_0$ such that $\left[H^1(\CP^2, E(-1))\right]
= \bv$ and $\operatorname{Supp} E_0/E = C$. Equivalently, it
parametrises $\Gamma$-equivariant $0$-dimensional quotients $E_0
\to Q$ such that $\left[H^0(\CP^2, Q)\right] = \bv - \bv^0$ and
$\operatorname{Supp}Q = C$. Such quotients depend only on a local
structure on $E_0$, so we can replace $E_0$ by
$W_s\otimes_\C{\mathcal O}_{\CP^2}$, where $W_s$ is the fiber of
$E_0$ at the origin considered as a representation of $\Gamma$.
The isomorphism class $\bw_s$ of $W_s$ is given by $\bw_s = \bw -
\bC\bv^0$, where $\bC$ is the class of the virtural representation
$\Wedge^0 Q - \Wedge^1 Q + \Wedge^2 Q = 2\rho_0 - Q$, and
$\bC\bv^0$ means the tensor product $\bC\otimes \bv^0$. Therefore
it becomes clear now that we have
\begin{equation*}
   \M(\bv,\bw)_x \cong
   \La(\bv-\bv^0 - m\delta, \bw_s) \times
   \prod_i \operatorname{Hilb}_0^{\lambda_i}{(\C^2)},
\end{equation*}
where $\delta$ is considered as the class of the regular
representation, and $\operatorname{Hilb}_0^{\lambda_i}{(\C^2)}$ is the
punctural Hilbert scheme, i.e., the inverse image of $\lambda_i [0]$
by the Hilbert-Chow morphism $\pi\colon \Hilb{(\C^2)}{\lambda_i}\to
S^{\lambda_i}(\C^2)$. The punctural Hilbert schemes are known to be
irreducible, thus $\M(\bv,\bw)_x$ is pure-dimensional if and only if
$\La(\bv-\bv^0-m\delta,\bw_s)$ is so. But the latter statement is
known \cite[\S5]{Na:1994}.

\section{A geometric construction of the affine Lie algebra}\label{sec:affine}
\setzero\vskip-5mm \hspace{5mm }

After writing \cite{KN}, the author tried to use this generalized ADHM
description to study these varieties $\M(\bv,\bw)$. But it turned out
to be not so easy as he had originally hoped.
When he struggled the problem, he heard a talk by Lusztig in ICM 90
Kyoto on a construction of canonical bases by using quivers. Lusztig's
construction \cite{Lu:can} was motivated by Ringel's construction
\cite{Ri} of the upper half part of the quantized enveloping algebra
via the Hall algebra.
The author thought that this construction should be useful to
attack the problem. Two years later, he began to understand the
picture. Quiver varieties $\M(\bv,\bw)$ are, very {\it roughly\/},
cotangent bundles of varieties used by Ringel and Lusztig, and
similar construction is possible \cite{Na:1994}. A little later,
he graduately realized that quiver varieties are also similar to
cotangent bundles of flag varieties and the map $\pi$ is an analog
of Springer resolution. These varieties had been used to give
geometric constructions of Weyl groups (Springer representations)
and affine Hecke algebras (Deligne-Langlands conjecture). (See a
beautifully written text book by N.~Chriss and V.~Ginzburg
\cite{CG} and the references therein for these matereial.)
The technique is the convolution product (see below) and works quite
general. So he (and some others) conjectured that these construction
should be adapted to quiver varieties. This conjecture turned out to
be true \cite{Na:1998,Na-qaff}. We explain the constructions in this
and next sections. The relation between our constructions and
Ringel-Lusztig construction was explained in \cite{Na:suugaku} and
will not be reproduced here.

\subsection{Convolution algebra}\vskip-5mm \hspace{5mm }

We apply the theory of the convolution algebra to varieties introduced
in the previous sections to obtain the universal enveloping algebra
$\bU(\ag)$ of the affine algebra $\ag$.

We continue to fix a representation $W$ of $\Gamma$ and denote by
$\bw$ its isomorphism class. For the $\Gamma$-fixed point set of
Hilbert schemes, studied in \secref{sec:Hilb}, $W$ is the trivial
representation.

We introduce the following notation:
\begin{equation*}
\begin{gathered}
   \M(\bw) \defeq \bigsqcup_n \M(n,r)^\Gamma
   = \bigsqcup_\bv \M(\bv,\bw),
\qquad
   \La(\bw) \defeq \bigsqcup_\bv \La(\bv,\bw),
\\
   \M_0(\infty,\bw) \defeq \bigcup_n \M_0(n,r)^\Gamma.
\end{gathered}
\end{equation*}
The first and second are disjoint union.
For the last, we use the inclusion
$\M_0(n,r)^\Gamma\subset \M_0(n',r)^\Gamma$ for $n\le n'$ given by
\begin{equation*}
   (E,\varphi,C) \mapsto \left(E,\varphi,C+(n'-n) 0\right).
\end{equation*}
For $r = 1$, these are
\begin{equation*}
   \bigsqcup_n \left(\Hilbn{(\C^2)}\right)^\Gamma, \quad
   \bigcup_n (S^n \C^2)^\Gamma = \bigcup_n S^n (\C^2/\Gamma),
\end{equation*}
where the inclusion $S^n (\C^2/\Gamma)\subset S^{n'}(\C^2/\Gamma)$ is
given by adding $(n'-n)0$ as above.

Rigorously speaking, we cannot study $\M(\bw)$ and $\M_0(\infty,\bw)$
directly since they are infinite dimensional. We need to work
individual spaces $\M(\bv,\bw)$, $\M_0(n,r)^\Gamma$. But we use those
spaces as if they are finite dimensional spaces for a notational
convenience.

We consider the fiber product
\begin{equation*}
   \Zw \defeq \M(\bw)\times_{\M_0(\infty,\bw)} \M(\bw).
\end{equation*}
It consists of pairs $(E,\varphi)$, $(E',\varphi')$ such that
\begin{enumerate}
\item $E^{\vee\vee} \cong E^{\prime\vee\vee}$
\item $\operatorname{Supp} E^{\vee\vee}$ and $\operatorname{Supp}
E^{\prime\vee\vee}$ are equal in the complement of the origin.
\end{enumerate}
The multiplicities of $\operatorname{Supp} E^{\vee\vee}$ and
$\operatorname{Supp} E^{\prime\vee\vee}$ at the origin may be
different since we consider the inclusion above.

One can show that this is a lagrangian subvariety in $\M(\bw)\times
\M(\bw)$. (The same remark as $\M(\bv,\bw)_x$ in \secref{sec:stratum}
applies here also.) Let us consider its top degree Borel-Moore
homology group
\begin{equation*}
    H_{\topdeg}(\Zw,\C).
\end{equation*}
More precisely, it is the subspace of
\begin{equation*}
   \prod_{\bv^1, \bv^2} H_{\topdeg}
   (\Zw\cap\left(\M(\bv^1,\bw)\times\M(\bv^2,\bw)\right),\C)
\end{equation*}
consisting of elements $(F_{\bv^1,\bv^2})$ such that
\begin{enumerate}
\item for fixed $\bv^1$, $F_{\bv^1,\bv^2} = 0$ for all but finitely
many choices of $\bv^2$,
\item for fixed $\bv^2$, $F_{\bv^1,\bv^2} = 0$  for all but finitely
many choices of $\bv^1$.
\end{enumerate}
The degree $\topdeg$ depends on $\bv^1$, $\bv^2$, but we supress the
dependency for brevity.

Let us consider the convolution product
\begin{equation*}
   \ast \colon H_{\topdeg}(\Zw,\C)\otimes H_{\topdeg}(\Zw,\C)
   \to H_{\topdeg}(\Zw,\C)
\end{equation*}
given by
\begin{equation*}
   c\ast c' = p_{13*}\left(p_{12}^*(c)\cap p_{23}^*(c')\right),
\end{equation*}
where $p_{ij}$ is the projection from the triple product
$\M(\bw)\times\M(\bw)\times\M(\bw)$ to the double product
$\M(\bw)\times\M(\bw)$ of the $i$th and $j$th factors. More detail for
the definition of the convolution product, say $p_{12}^*$, $\cap$, is
explained in \cite{CG}, but we want to emphasize one point. The
statement that the result $c\ast c'$ has top degree is the consequence
of $\dim Z(\bw) = \frac12 \dim \M(\bw)\times\M(\bw)$.
Although we are considering $\Zw$ having infinitely many connected
components, the convolution is well-defined and $H_{\topdeg}(\Zw,\C)$
is an associative algebra with unit, thanks to the above definition of
the subspace in the direct product.

For $x\in \M_0(\infty, \bw)$, let $\M(\bw)_x$ be the inverse image
$\pi^{-1}(x)$ in $\M(\bw)$. We consider the top degree homology group
\begin{equation*}
   H_{\topdeg}(\M(\bw)_x, \C)
\end{equation*}
which is the usual direct sum of
$H_{\topdeg}(\M(\bw)_x\cap\M(\bv,\bw),\C)$ (unlike the case of $\Zw$).
The convolution product makes this space into a module of
$H_{\topdeg}(\Zw,\C)$.

\begin{Theorem}
Let $\bU(\ag)$ be the universal enveloping algebra of the untwisted
affine algebra $\ag$ corresponding to $\Gamma$.
\textup({\bf NB}\textup: not a `quantum' version\textup). There exists
an algebra homomorphism
\begin{equation*}
   \bU(\ag) \to H_{\topdeg}(\Zw,\C).
\end{equation*}
Furthermore, if we consider $H_{\topdeg}(\M(\bw)_x, \C)$ as a
$\bU(\ag)$-module via the homomorphism, it is an irreducible
integrable highest weight representation and the direct summands
$H_{\topdeg}(\M(\bw)_x\cap\M(\bv,\bw),\C)$ are weight spaces.
\end{Theorem}

This theorem was essentially proved in \cite{Na:1994} with a
modification for general $x$ mentioned above.

The highest weight of $H_{\topdeg}(\M(\bw)_x, \C)$ and weights of
$H_{\topdeg}(\M(\bw)_x\cap\M(\bv,\bw),\C)$ are determined explicitly
in terms of $\bv$, $\bw$ and the stratum to which $x$ belongs.
For example, in the case of the central fiber $\La(\bv,\bw)$, the
highest weight is $\bw$, considered as a dominant integral weight as
\(
   \bw = \sum_i w_i \Lambda_i,
\)
where $w_i$ is the $\rho_i$ component of $\bw$, and $\Lambda_i$ is the
$i$th fundamental weight. Here we use the identification of the
irreducible representation $\rho_i$ and a vertex of the affine Dynkin
diagram given by McKay correspondence. The weight of
$H_{\topdeg}(\La(\bv,\bw),\C)$ is $\bw - \bv$, where
\(
   \bv = \sum_i v_i \alpha_i
\)
with the $\rho_i$-component $v_i$ of $\bv$ and $i$th simple root
$\alpha_i$.
The highest weight vector is the fundamental class $[\La(0,\bw)]$,
where $\M(0,\bw) = \La(0,\bw)$ consists of a single point
$E = W\otimes_\C{\mathcal O}_{\CP^2}$.

For the case studied in \secref{sec:Hilb}, $\bw$ is the $0$th
fundamental weight $\Lambda_0$. The corresponding integrable highest
weight representation is called the {\it basic\/} representation in
literature. If we vary $\bw$, we get all integrable highest weight
representations as $\bigoplus_\bv H_{\topdeg}(\La(\bv,\bw),\C)$. It is
worth while remarking that this is an extension of \eqref{eq:Cartan}
since the Cartan subalgebra $\mathfrak h$ is naturally contained in
the Cartan subalgebra. Furthermore, the finite dimensional Lie algebra
$\g$ is embedded in the basic representation, and we get
\begin{equation*}
   \g \cong \bigoplus_{v_0 = 1} H_{\topdeg}(M(\bv),\C),
\end{equation*}
where $v_0$ is the $\rho_0$-isotropic component of $\bv$. This is an
extension of \eqref{eq:Cartan}, mentioned before. In fact, it is easy
to see that if $\bv$ is not $\delta$, then $M(\bv)$ is either empty,
or a single point. The latter holds if and only if $\bv$, considered
as an element of $\mathfrak h^*$ by removing $v_0$, is a root of $\g$.

If we fix $\bw$ and vary $x$, we still obtain various integrable
highest weight representations. The highest weight of
$H_{\topdeg}(\M(\bw)_x,\C)$ is $\bw - \bv^0 - m\delta$, where $\bv^0$,
$m$ are determined by $x$ as in \secref{sec:stratum}, and
$\bw$, $\bv^0$ are considered as weights as above. The weight of
$H_{\topdeg}(\M(\bw)_x\cap\M(\bv,\bw),\C)$ is equal to $\bw -
\bv$.
All of their highest weights are less than or equal to $\bw$ with
respect to the dominance order. In particular, when $\bw = \Lambda_0$,
those have highest weights $\Lambda_0 - n\delta$ for some $n\in\Z_{\ge
0}$. They are essentially isomorphic to the basic representation.

We explain how the algebra homomorphism \( \bU(\ag) \to
H_{\topdeg}(\Zw,\C) \) is defined. It is enough to define the image of
Chevalley generators $e_i$, $f_i$, $h_i$ ($i\in I$), $d$ of $\bU(\ag)$
(and check the defining relations). The images of $h_i$ and $d$ are
multiples of fundamental classes of diagonales in
$\M(\bw)\times\M(\bw)$. More precisely, the multiple is determined so
that the weight of $H_{\topdeg}(\M(\bw)_x\cap\M(\bv,\bw),\C)$ is equal
to $\bw - \bv$.
The image of $e_i$ is the fundamental classe of the so-called
`Hecke correspondence': \setcounter{equation}{0}
\begin{equation}\label{eq:Hecke}
   \bigsqcup_\bv \left\{\left. \left((E,\varphi), (E', \varphi')\right)
   \in \M(\bv,\bw)\times\M(\bv+\rho_i,\bw)\, \right|
   E\subset E' \right\}.
\end{equation}
It is known that each component is a nonsingular lagrangian subvariety of
$\M(\bv,\bw)\times\M(\bv+\rho_i,\bw)$. Hence it is an irreducible
component of $\Zw$. The image of $f_i$ is given by swapping the first
and second factors, up to sign.

As an application of the above construction, we get a base of
$H_{\topdeg}(\M(\bw)_x, \C)$ indexed by the irreducible components of
$\M(\bw)_x$. It has a structure of the crystal in the sense of
Kashiwara, and is isomorphic to the crystal of the corresponding
integrable highest weight module of the quantum affine algebra by
Kashiwara-Saito \cite{KS,Saito}. (See also \cite{Na:tensor} for a
different proof.) However, the base itself is different from the
specialization of the canonical (= global crystal) base of the quantum
affine algebra module at $q = 1$. A counter example was found in
\cite{KS}. The base given by irreducible components is named {\it
semicanonical base\/} by Lusztig \cite{Lu:semican}.

\subsection{Lower degree homology groups}\vskip-5mm \hspace{5mm }

The construction in the previous subsection, in fact, gives us a
structure of a representation of $\ag$ on
\[
   H_{\topdeg - d} (\M(\bw)_x, \C)
\]
for each fixed integer $d$. It is an integrable representation, and
decomposes into irreducible representations.
The multiplicity formula can be expressed in terms of the intersection
cohomology thanks to Beilinson-Bernstein-Deligne-Gabber's
decomposition theorem~\cite{BBD} applied to the morphism $\pi\colon
\M(\bw)\to \M_0(\infty,\bw)$. In our situation, $\pi$ is a {\it
semi-small\/} morphism, i.e., the restriction of $\pi$ to the inverse
image of the stratum \eqref{eq:strat} is a topological fiber bundle,
and
\[
   2\dim \M(\bw)_x \le \codim \mathcal O_x,
\]
where $\mathcal O_x$ is the stratum containing $x$. Then as
observed by Borho-MacPherson~\cite{BM}, the decomposition theorem
is simplified. We introduce several notation to describe the
formula. We choose a point $y$ from each stratum in
\eqref{eq:strat}. We denote the stratum containing $y$ by
$\mathcal O_y$. Let $IC(\mathcal O_y)$ is the intersection
homology complex of $\mathcal O_y$ with respect to the trivial
local system. \setcounter{equation}{1}
\begin{Theorem}\label{thm:decomp}
We have the following decomposition as a representation of $\ag$\textup:
\begin{equation*}
   H_{\topdeg - d} (\M(\bw)_x, \C)
   = \bigoplus_{y} H^{d + \dim {\mathcal O}_x}(i_x^! IC(\mathcal O_y))
     \otimes H_{\topdeg}(\M(\bw)_y, \C),
\end{equation*}
where $i_x\colon \{x\}\to \M_0(\infty,\bw)$ is the inclusion. Here
$\ag$ acts trivially on the first factor of the right hand side.
\end{Theorem}

For a general semi-small morphism, we may have an intersection
homology complex with respect to a {\it nontrivial\/} local system
in the decomposition. In order to show that such a summand does
not appear, the fact that $H_{\topdeg}(\M(\bw)_y, \C)$ is a
highest weight module plays a crucial role (see \cite{Na-qaff} for
detail).

Note also that when $\bw = \rho_0$, the closure of each stratum is a
symmetric product of $\C^2/\Gamma$ \eqref{eq:strat}. In particular,
they only have finite quotient singularities, and their intersection
homology complex are equal to the constant sheaf. Therefore our
formula are simplified. One finds that $H_*(\La(\bw),\C)$ is
isomorphic to the so-called `Fock space'.
Later we will show that the total homology $H_*(\La(\bw),\C)$ has
a structure of a representation of a specialized quantum toroidal
algebra in the next section. Then this observation generalizes a
result \cite{VV:tor,STU} for type $A_n^{(1)}$ to untwisted affine
Lie algebras of type $ADE$. \pagebreak
\section{Equivariant {\boldmath $K$}-theory and quantum toroidal algebras}
\setzero\vskip-5mm \hspace{5mm }

In this section, we replace the top degree homology group
$H_{\topdeg}$ in the previous section by equivariant $K$-groups. Then
we obtain a geometric construction of a quantum toroidal algebra
$\Ut$. It is a $q$--analog of the loop algebra ${\mathbf L}{\ag} =
\ag\otimes_\C \C[z, z^{-1}]$ of the affine algebra $\ag$. Since $\ag$
is already a (central extension of) loop alebra of $\g$, ${\mathbf
L}{\ag}$ is a `double-loop' algebra of $\g$.
The quantum toroidal algebra is defined by replacing $\g$ by $\ag$
in the so-called Drinfeld realization of the quantum loop algebra
$\Ul$, which is a subquotient of the quantum affine algebra $\Ua$,
defined by Drinfeld, Jimbo.

Let $G_\bw = \operatorname{Aut}_\Gamma(W)$ be the group of
automorphisms of the $\Gamma$-module $W$. If $w_i$ is the multiplicity
of $\rho_i$ in $W$, we have $G_\bw \cong \prod_i \GL_{w_i}(\C)$. We
have a natural action of $G_\bw$ on $\M(\bw)$ and $\M_0(\infty,\bw)$
by the change of the framing:
\begin{equation*}
   \varphi \mapsto g\circ\varphi, \qquad g\in G_\bw.
\end{equation*}
The projective morphism $\M(\bw)\to\M_0(\infty,\bw)$ is equivariant.

Let $\C^*$ act on $\C^2$ by $t\cdot(x,y) = (tx, ty)$. It extends to an
action on $\CP^2$, where it acts trivially on $\linf$. Note that this
action commutes with the $\Gamma$-action. Then we have a natural
induced $\C^*$-action on $\M(\bw)$ and $\M_0(\infty,\bw)$ so that the
projection $\pi$ is equivariant. Combining two actions we have an
action of $G_\bw\times\C^*$ on $\M(\bw)$ and $\M_0(\infty,\bw)$.
(This action is different from the action studied in \cite{Na-qaff}
for type $A_1^{(1)}$. We need to change the definition of $\Ut$ in
that paper to apply the result in this section. This comes from the
umbiguity of the definition of a $q$-analog of the Cartan matrix.)

Let $K^{G_\bw\times\C^*}(\Zw)$ be the equivariant $K$-homology
group of $\Zw$ with respect to the above $G_\bw\times\C^*$-action.
(More precisely, it should be defined as a subspace of the direct
product as in the case of homology groups.) It is a module over
the representation ring $R(G_\bw\times\C^*) =
\Z[q,q^{-1}]\otimes_\Z R(G_\bw)$, where $q$ is the natural
$1$-dimensional representation of $\C^*$. The convolution product
makes $K^{G_\bw\times\C^*}(\Zw)$ into a
$R(G_\bw\times\C^*)$-algebra. We divide its torsion part over
$\Z[q,q^{-1}]$ and denote it by
$K^{G_\bw\times\C^*}(\Zw)/\operatorname{torsion}$. (It is
conjectured that the torsion is, in fact, $0$.)

\begin{Theorem}
There exists a $\Z[q,q^{-1}]$-algebra homomorphism
\begin{equation*}
   \Uti \to K^{G_\bw\times\C^*}(\Zw)/\operatorname{torsion},
\end{equation*}
where $\Uti$ is a certain $\Z[q,q^{-1}]$-subalgebra
\textup(conjecturally an integral form\textup) of $\Ut$.
\end{Theorem}

The definition of the homomorphism is similar to the case of homology
groups. The image of the $q$-analog of $e_i\otimes z^r$ is given by a
natural line bundles on the Hecke correspondence \eqref{eq:Hecke}
whose fiber at $((E,\varphi), (E',\varphi'))$ is $H^0(E'/E)^{\otimes r}$.

Let us explain how we can use this algebra homomorphism to study
representations of specialized quantum toroidal algebra
\(
   \Ute = \Uti|_{q=\ve}
\),
where $\ve$ is a nonzero complex number which may or may not be a root of
unity. A natural generalization of finite dimensional representations
of $\Ule$ are {\it l\/}--integrable representations. (See
\cite{Na-qaff} for the definition.) The Drinfeld-Chari-Pressley
classification \cite{Dr,CP} of irreducible finite dimensional
reprentation of $\Ule$ has a natural analog in $\Ute$: Irreducible
{\it l\/}--integrable representations of $\Ute$ are parametrized by
$I$--tuple of polynomials $P_i(u)$ with $P_i(0) = 1$, where $I$ is
the set of vertices of the affine Dynkin diagram.

Irreducible representations are obtained in the following way.
Let us consider the equivariant homology
$K^{G_\bw\times\C^*}(\La(\bw))$, which is without torsion. It is a
module of $\Uti$ and called a {\it universal standard module}. For a
semisimple element $(s,\ve)\in G_\bw\times\C^*$, we consider the
evaluation homomorphism $R(G_\bw\times\C^*)\to \C$. Then the
specialization
\[
   K^{G_\bw\times\C^*}(\La(\bw))\otimes_{R(G_\bw\times\C^*)}\C
\]
is a representation of $\Ute$. This is called a {\it standard module}.
It has a unique irreducible quotient, and the associated polynomials
are the characteristic polynomials of components of $s$. (Recall
$G_\bw = \prod_{i\in I} \GL_{w_i}(\C)$.)

In order to state character formulas, which is very similar to
\thmref{thm:decomp}, we need a little more notation. Let $A$ be the
Zariski closure of powers of $(s,\ve)$ in $G_\bw\times\C^*$. Let
$\M(\bw)^A$, $\M_0(\infty,\bw)^A$, $\Zw^A$ be the fixed point sets.
We have a chain of natural algebra homomorphisms
\begin{multline*}
   K^{G_\bw\times\C^*}(\Zw)_{R(G_\bw\times\C^*)}\C
        \to
   K^{A}(\Zw)_{R(A)}\C
\\
        \to
   K(\Zw^A)\otimes_\Z \C
        \to
   H_*(\Zw^A,\C),
\end{multline*}
where the first one is induced by the inclusion $A\subset G_\bw\times\C^*$,
the second one is given by the localization theorem in the equivariant
$K$-theory, and the last one is the Chern character homomorphism.
(In fact, we need `twists' for the last two. See \cite{CG} for detail.)

There exists a natural stratification of
\(
   \M_0(\infty, \bw)^A
\)
similar to \eqref{eq:strat}. We choose a point $y$ in each stratum and
denote by $\mathcal O_y$ the stratum containing $y$.
If $\M(\bw)^A_y$ denotes the inverse image of $y$ in $\M(\bw)^A$
under $\pi$, the homology group $H_*(\M(\bw)^A_y,\C)$ is a
representation of $H_*(\Zw^A,\C)$, and hence that of $\Ute$. (When
$y=0$, it is the standard module.) Analog of \thmref{thm:decomp}
is the following:
\begin{Theorem}\label{thm:mult}
We have the following in the Grothendieck group of the abelian
category of {\it l\/}--integrable representations of $\Ute$
\begin{equation*}
   H_*(\M(\bw)^A_x,\C)
   = \sum_y H^*(i_x^! IC(\mathcal O_y)) \otimes L_y,
\end{equation*}
where $L_y$ is the unique irreducible quotient of $H_*(\M(\bw)^A_y,
\C)$.
\end{Theorem}
\noindent
(The right hand side is an infinite sum, so we
must understand it in an appropriate way. But it should be clear how
it can be done.)

In this case, we apply the decomposition theorem to $\M(\bw)^A\to
\M_0(\infty,\bw)^A$. It is {\it not\/} semi-small any more. So the
degrees in the left and righ hand sides do not have clear relations.

\begin{Remarks}
(1) We stated our results for the affine Lie algebra $\ag$ and the
quantum toroidal algebra $\Ut$. But they also hold for the finite
dimensional Lie algebra $\ag$ and the quantum loop algebra $\Ul$, if
we impose the condition $v_0 = w_0 = 0$. It is known that $\Ul$ is a
Hopf algebra (since it is a subquotient of the quantum affine
algebra), and the standard modules are isomorphic to tensor products
of {\it l\/}--fundamental representations when $\ve$ is not a root of
unity \cite{VV:std}. Here {\it l\/}--fundamental representations are
irreducible representations corresponding to $\bw =
\rho_i$. In particular, the tensor product decomposition in the
representation ring can be expressed in terms of intersection homology
groups.

(2) Our remaining tasks are computing dimensions of
$H^*(i_x^!IC({\mathcal O}_y))$ appearing in
Theorems~\ref{thm:decomp}, \ref{thm:mult}. In
\cite{Na:qchar,Na:qchar-main} we gave a purely combinatorial
algorithm to compute them. This algorithm can be made into a
computer program.
The algorithm was stated for the quantum loop algebra $\Ul$, but works
also for $\Ut$. This means that for any given stratum $\mathcal O_y$,
$\dim H^*(i_x^!IC({\mathcal O}_y))$ is, in principle,
computable. However, it is practically difficult to compute because we
need lots of memory. And, for $\Ut$, the summation is infinite. So
having an algorithm to compute each term is not a strong statement. It
is desirable to have an alternative method to compute them. That has
be done in some special classes of representations \cite{Na:KR}.
\end{Remarks}

\providecommand{\bysame}{\leavevmode\hbox to3em{\hrulefill}\thinspace}

\label{lastpage}

\end{document}